\documentclass[a4paper,12pt]{elsarticle}
\usepackage{amsmath,amsfonts,amssymb,mathtools}
\usepackage{hyperref}
\usepackage{numcompress}
\makeatletter
\def\ps@pprintTitle{%
   \let\@oddhead\@empty
   \let\@evenhead\@empty
   \let\@oddfoot\@empty
   \let\@evenfoot\@oddfoot
}
\makeatother
\biboptions{square,comma}
\begin{document}
\begin{frontmatter}
\title{A comparison between numerical solutions to fractional differential equations: Adams-type predictor-corrector and multi-step generalized differential transform method}
\author{Alireza Momenzadeh}
\author{Sima Sarv Ahrabi\corref{cor1}}
\cortext[cor1]{Corresponding author}
\begin{abstract}
In this note, two numerical methods of solving fractional differential equations (FDEs) are briefly described, namely predictor-corrector approach of Adams-Bashforth-Moulton type and multi-step generalized differential transform method (MSGDTM), and then a demonstrating example is given to compare the results of the methods. It is shown that the MSGDTM, which is an enhancement of the generalized differential transform method, neglects the effect of non-local structure of fractional differentiation operators and fails to accurately solve the FDEs over large domains.
\end{abstract}


\end{frontmatter}

\section{Introduction}
Fractional calculus has been studied by mathematicians for years and fractional differential equations (FDEs) specifically proven crucial to the more accurate mathematical modeling of many physical phenomena in diverse branches of science. The nature of most physical systems is governed by after-effect or memory, in other words they depend on the all previous time history instead of just for an instant and could be elegantly modeled by using FDEs. Frequent utilization of fractional differential equations could be observed in, for instance, the theory of viscoplasticity \cite{diethelm1998fracpece,carpinteri2014fractals}, physics \cite{atangana2016new,nasholm2013fractional}, civil engineering \cite{lopes2013fractional,koh1990application}, mechanics of solids \cite{rossikhin1997applications,rossikhin2010application}, control theory \cite{agrawal2004general,jelicic2009optimality}, biological systems \cite{magin2010fractional,ding2009fractional,glockle1995fractional,magin2006fractional},etc. However, such differential equations of fractional order, specifically those possessing nonlinear terms might be impossible to be solved analytically. Consequently, numerical methods are mainly developed from classical approaches, which are appropriately modified to be able to deal with the fractional differential operators. Method of extrapolation \cite{diethelm1997numerical,chen2008adi}, fractional linear multi step method \cite{lubich1985fractional,galeone2006multistep,galeone2009explicit}, quadrature-based direct method \cite{baleanu2016fractional}, predictor-corrector approach of Adams-Bashforth-Moulton type for numerical solution of FDEs \cite{diethelm1998fracpece,baleanu2016fractional,diethelm1997algorithm,diethelm2004detailed,diethelm2003efficient,garrappa2011predictor,diethelm2010analysis} and fractional differential transform method \cite{arikoglu2007solution,odibat2008generalized} provide numerical solutions for linear and nonlinear  differential equations of fractional order.
In sections \ref{predictor} and \ref{msgdtm}, a brief description of predictor-corrector method of Adams-Bashforth-Moulton and multi-step fractional differential transform method (MSGDTM) are presented, then in section \ref{compare}, an illustrating example is given, in which a fractional differential equation is solved by using the two methods and the results will be graphically compared with each other.
\section{Predictor-Corrector approach of Adams type}
\label{predictor}
The algorithm that will be briefly described could be considered as a fractional variant of the classical second-order Adams-Bashforth-Moulton method \cite{diethelm1998fracpece,diethelm2010analysis}. In this note, the main emphasis will be placed on the single-term Caputo fractional differential equations \cite{diethelm2010analysis} for $0<\alpha \le 1$ , where $\alpha $  is the order of the fractional derivative.
Consider the initial value problem of fractional differential equation:
\begin{equation}
\label{ini1}
\begin{cases}
D_{{t}_{0}}^{\alpha} y\left (t \right)=f\left (t,y\left( t \right) \right),\\
y\left( {{t}_{0}} \right)={{y}_{0}}.
\end{cases}
\end{equation}
In order to assure the existence and uniqueness of the solution to the problem (\ref{ini1}), it is assumed that $f\left( t,y \right)$ is continuous and fulfils a Lipschitz condition with respect to the second variable \cite[Theorem~6.5]{diethelm2010analysis}. The initial value problem (\ref{ini1}) could be converted into an equivalent Volterra integral equation \cite{kilbas2006theory}
\begin{equation}
\label{volterra}
y\left( t \right)={{y}_{0}}+\frac{1}{\Gamma \left( \alpha  \right)}{\int_{{t_0}}^{t}(t-s)^{\alpha -1}}f\left( s,y\left( s \right) \right)ds.\
\end{equation}
The method presents a numerical approach in solving Eq.(\ref{volterra}) and is said to be PECE (Predict, Evaluate, Correct and Evaluate) type because an initial approximation $y_{k}^{P}$ , the so-called predictor, is first evaluated for the solution
\begin{equation}
\label{voltraP}
y_{k}^{P}={{y}_{0}}+\frac{1}{\Gamma \left( \alpha  \right)}\sum\limits_{j=0}^{k-1}{{{b}_{j,k}}f\left( {{t}_{j}},{{y}_{j}} \right)},\
\end{equation}
where $k\in {{\mathbb{Z}}^{+}}$ and the nodes ${{t}_{j}}\text{ }\left( j=0,1,...,k \right)$ are used to calculate $y_{k}^{P}$  and consequently ${{y}_{k}}$ step by step. Considering the equal step size with some fixed $h$ (${{t}_{j}}=jh$), the weight ${{b}_{j,k}}$ is evaluated by the formula below,
\begin{equation}
\label{bj}
{{b}_{j,k}}=\frac{{{h}^{\alpha }}}{\alpha }\left( {{\left( k-j \right)}^{\alpha }}-{{\left( k-1-j \right)}^{\alpha }} \right).
\end{equation}
Then the method gives the corrector formula, which is
\begin{equation}
\label{yk}
{{y}_{k}}={{y}_{0}}+\frac{1}{\Gamma \left( \alpha  \right)}\left( {{a}_{k,k}}f\left( {{t}_{k}},y_{k}^{P} \right)+\sum\limits_{j=0}^{k-1}{{{a}_{j,k}}f\left( {{t}_{j}},{{y}_{j}} \right)} \right),
\end{equation}
where the weight ${{a}_{j,k}}$ is given by
\begin{alignat}{2}
\label{akk}
\begin{dcases}
   {{a}_{0,k}}=\frac{{{h}^{\alpha }}}{\alpha \left( \alpha +1 \right)}\left( {{\left( k-1 \right)}^{\alpha +1}}-{{k}^{\alpha }}\left( k-1-\alpha  \right) \right),\\
   \begin{aligned}
    {{a}_{j,k}}=\frac{{{h}^{\alpha }}}{\alpha \left( \alpha +1 \right)}\left( {{\left( k+1-j \right)}^{\alpha +1}}+{{\left( k-1-j \right)}^{\alpha +1}}-2{{\left( k-j \right)}^{\alpha +1}} \right),\\
    1\le j\le k-1,\\
    \end{aligned}\\
        {{a}_{k,k}}=\frac{{{h}^{\alpha }}}{\alpha \left( \alpha +1 \right)}.
\end{dcases}
\end{alignat}
The basic algorithm, the fractional Adams-Bashforth-Moulton method, could be completely described by (\ref{voltraP}) and (\ref{yk}) with the weights ${{a}_{j,k}}$ and ${{b}_{j,k}}$ defined according to (\ref{bj}) and (\ref{akk}).
\section{Multi-step generalized differential transform method}
\label{msgdtm}
The differential transform method (DTM) has been extended by \cite{arikoglu2007solution} in order to be able to obtain the approximate solutions of FDEs. The new technique is called as fractional differential transform method (FDTM). Then, \cite{odibat2008generalized} presented a new generalization of one-dimensional differential transform method in order to deal with differential equations of fractional orders. The method is mentioned as generalized differential transform method (GDTM), which is based on generalized Taylor's formula and Caputo fractional derivative. For convenience, the method is briefly described below. The detailed description of the technique could be observed in the main articles \cite{arikoglu2007solution,odibat2008generalized}. The GDTM results in representing the solution to Eq. (\ref{ini1}) as the form below:
\begin{equation}
\label{8}
y\left( t \right)=\sum\limits_{k=0}^{\infty}{Y\left( k \right){{\left( t-{{t}_{0}} \right)}^{k\alpha }}}.
\end{equation}
The coefficients $Y\left(k\right)$ could be straightforwardly evaluated by using the recurrence equation
\begin{equation}
Y\left(k+1\right)=\frac{\Gamma \left( \alpha k+1 \right)}{\Gamma \left( \alpha \left( k+1 \right)+1 \right)}F\left( k,Y\left( k \right) \right),
\end{equation}
where $Y\left( 0 \right)$ is assessed to be equal to $y\left( {{t}_{0}} \right)$ , i.e. $Y\left( 0 \right)={{y}_{0}}$ and the term $F\left( k,Y\left( k \right) \right)$, which is mentioned as differential transform of $f\left (t,y\left( t \right) \right)$, could be in succession determined by using the techniques, provided in related articles \cite{arikoglu2007solution,odibat2008generalized,odibat2008differential,bervillier2012status}.
The important aspect of this method is the easy generalization of the Taylor method by providing recurrence schemes for problems involving fractional derivatives.
Nonetheless, the method is subjected to the wrong assumption that the Eq. (\ref{8}) is sufficient to provide a suitable approximation of the exact solution of the FDEs. It has been definitely proved that, in the general circumstances, the exact solution of a differential equation of fractional order is expanded in mixed powers, i.e. a collection of integer and non-integer orders \cite{diethelm2010analysis}. In order to expand upon the point, it could be mentioned that the Volterra integral equation (\ref{volterra}) has been demonstrated to have an asymptotic expansion in terms of mixed powers of $\left( t-{{t}_{0}} \right)$ and ${{\left( t-{{t}_{0}} \right)}^{\alpha }}$ (see \cite{lubich1983runge}). Therefore, eliminating integer-order powers in the series expansion leads to a poor approximation.

Moreover,the method is originally based on Taylor series; therefore, the solution could be only relied on for a short period after the initial time. In order to deal with this restriction, in \cite{erturk2011approximate} the authors have developed multi-step generalized differential transform method (MSGDTM), which will be concisely described.
The multi-step generalized differential transform method is simply formed on the idea of dividing the time interval $[{{t}_{0}},T]$ into \emph{M} sub-intervals $[{{t}_{j}},{{t}_{j+1}}],\text{ }j=0,1,...,M-1$ of equal step size $h={\left( T-{{t}_{0}} \right)}/{M}\;$ and applying the generalized differential transform method (GDTM) to each of them. Over the first interval $[{{t}_{0}},{{t}_{1}}]$, the approximate solution ${{y}_{1}}\left( t \right)$ is obtained by using the initial condition ${{y}_{1}}\left( {{t}_{0}} \right)={{y}_{0}}$. In order to find the solution ${{y}_{2}}\left( t \right)$ over the second sub-interval $[{{t}_{1}},{{t}_{2}}]$  again, GDTM is applied to the Eq.(\ref{ini1}) with the initial condition ${{y}_{2}}\left( {{t}_{1}} \right)={{y}_{1}}\left( {{t}_{1}} \right)$, which is taken from the previous step. The same process is consistently repeated to the last sub-interval $[{{t}_{M-1}},{{t}_{M}}]$ and the MSGDTM assumes the \emph{N}th-order approximate solution to ultimately be
\begin{equation}
\ y(t)=\begin{cases}
       y_1(t)=y_0 +\sum\limits_{k=1}^{N}Y_1(k)(t-t_0)^{k\alpha},& t\in [t_0,t_1]\\
       y_2(t)=y_1(t_1)+\sum\limits_{k=1}^{N}Y_2(k)(t-t_1)^{k\alpha},& t\in [t_1,t_2]\\
       .\\
       .\\
       y_{M}(t)=y_{M-1}(t_{M-1})+\sum\limits_{k=1}^{N}Y_M(k)(t-t_{M-1})^{k\alpha},& t\in [t_{M-1},t_M].\\
        \end{cases}
\end{equation}
 This method is practically successful in solving ordinary differential equations (ODEs) due to the fact that the derivatives of integer orders are local in nature, i.e. the integer order derivative of a function $f\left( t \right)$ at a point $t$ depends only on the graph of $f\left( t \right)$  very close to the point $t$ . Unlike classical derivatives, fractional order differentiation has non-local property. Consider the ordinary differential equation
\begin{equation}
\label{ode}
\dot{y}=f\left( t,y \right),
\end{equation}
with the initial condition
\begin{equation}
\label{odein}
y({{t}_{0}})={{y}_{0}}.
\end{equation}
It is assumed that there exists a unique solution for the initial value problem (\ref{ode}) and (\ref{odein}). Suppose that ${{y}_{k}}$, an approximation of the solution at time ${{t}_{k}}$  has been already calculated. It is possible to evaluate ${{y}_{k+1}}$ , the approximate value of the solution at time ${{t}_{k+1}}$ , by using the equation
\begin{equation}
\label{k+1}
{{y}_{k+1}}={{y}_{k}}+\int_{{{t}_{k}}}^{{{t}_{k+1}}}{f\left( s,y\left( s \right) \right)ds}.
\end{equation}
The lower limit of integral in Eq.(\ref{k+1}) is the node ${{t}_{k}}$; therefore, it is possible to solve an ordinary differential equation at each step by using the initial condition, which comes from the previous step. Now, consider the fractional initial value problem (\ref{ini1}). The approximate value of the solution at time ${{t}_{k+1}}$ could be evaluated by
\begin{equation}
\label{gamma}
{{y}_{k+1}}={{y}_{0}}+\frac{1}{\Gamma \left( \alpha  \right)}{\int_{{{t}_{0}}}^{{{t}_{k+1}}}{\left( {{t}_{k+1}}-s \right)}^{\alpha -1}}f\left( s,y\left( s \right) \right)ds.
\end{equation}
It could be easily observed that the Eq.(\ref{gamma}) is different from Eq.(\ref{k+1}) and the interval of the integration in Eq.(\ref{gamma}) starts at the initial time ${{t}_{0}}$ instead of ${{t}_{k}}$. This is due to the non-local structure of the fractional differential operators. Now, it is obvious from the Eq.(\ref{gamma}) that dividing the interval $[{{t}_{0}},T]$ into several sub-intervals and solving the fractional differential equation for each one with initial conditions, which come from previous sub-intervals and neglecting the non-local property of fractional derivatives, results in solving a wrong equation for each sub-interval
\begin{equation}
{{y}_{k+1}}={{y}_{k}}+\frac{1}{\Gamma \left( \alpha  \right)}{\int_{{{t}_{k}}}^{{{t}_{k+1}}}{\left( {{t}_{k+1}}-s \right)}^{\alpha -1}}f\left( s,y\left( s \right) \right)ds,
\end{equation}
instead of the correct equation, i.e. Eq. (\ref{gamma}), which is rewritten in the form
\begin{multline}
\label{16}
{{y}_{k+1}}={{y}_{0}}+\frac{1}{\Gamma \left( \alpha  \right)}{\int_{{{t}_{0}}}^{{{t}_{k}}}{\left( {{t}_{k+1}}-s \right)}^{\alpha -1}}f\left( s,y\left( s \right) \right)ds+\frac{1}{\Gamma \left( \alpha  \right)}{\int_{{{t}_{k}}}^{{{t}_{k+1}}}{\left( {{t}_{k+1}}-s \right)}^{\alpha -1}}\\
\times f\left( s,y\left( s \right) \right)ds.
\end{multline}
It could be apparently perceived what will be neglected here is the term
\begin{equation}
\frac{1}{\Gamma \left( \alpha  \right)}{\int_{{{t}_{0}}}^{{{t}_{k}}}{\left( {{t}_{k+1}}-s \right)}^{\alpha -1}}f\left( s,y\left( s \right) \right)ds,
\end{equation}
in the right-hand side of the Eq.(\ref{16}), which consequently causes a discontinuity in the first derivative of the solution at the beginning of each sub-interval.
\section{Numerical results}
\label{compare}
\begin{figure}[ht!]
\centering
\includegraphics[scale=.7]{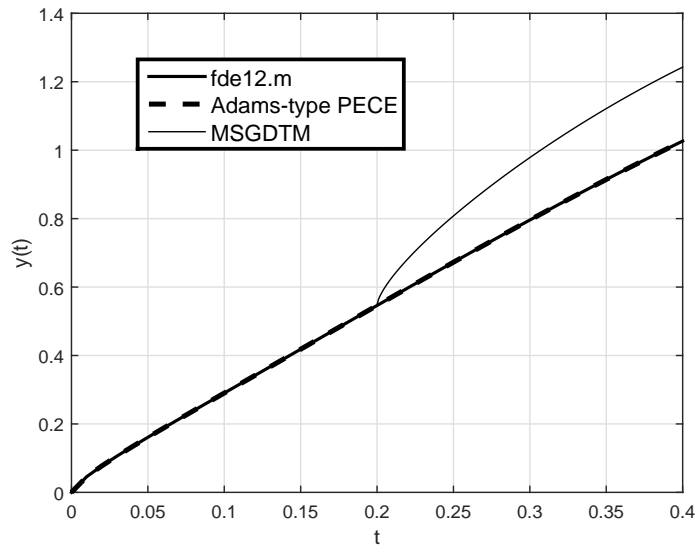}
\caption{Comparison of Adams-type predictor-corrector method and MSGDTM over two sub-intervals for $\alpha =0.7$}
\label{fig1}
\end{figure}

\begin{figure}[ht!]
\centering
\includegraphics[scale=.7]{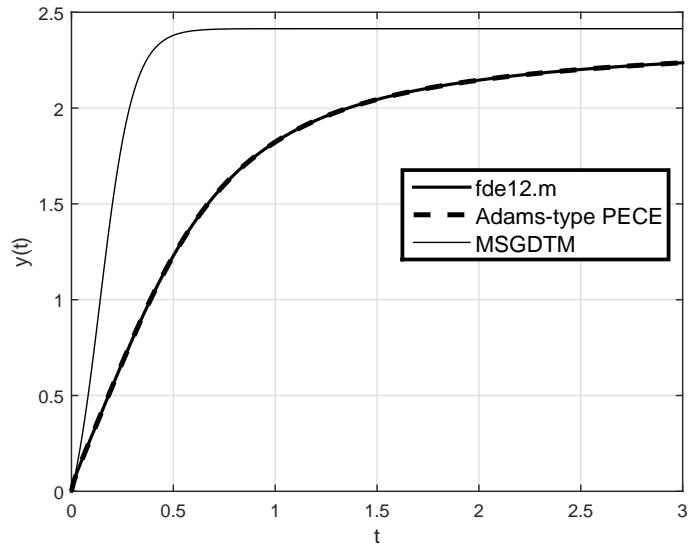}
\caption{Comparison of Adams-type predictor-corrector method $\left( h=0.01 \right)$ and MSGDTM $\left( M=300\right)$ for $\alpha =0.7$. }
\label{fig2}
\end{figure}

This section is allocated to an example demonstrating the numerical solution of a fractional differential equation by using the method of PECE of Adams-Bashforth-Moulton type and the MSGDTM.
Consider the fractional Riccati differential equation (see \cite{odibat2008generalized}):
\begin{equation}
\label{4.1}
   \begin{cases}
    D_{0}^{\alpha }y\left( t \right)=2y-{{y}^{2}}+1,&\text{        } t>0\text{   ,      } 0<\alpha \le 1,\\
    y\left( 0 \right)=0.
\end{cases}
\end{equation}
The objective is to follow the MSGDTM for the interval $I=[0,0.4]$ by dividing it into two sub-intervals, namely ${{I}_{1}}=[0,0.2]$ and ${{I}_{2}}=[0.2,0.4]$. The differential transformation of Eq. (\ref{4.1}) is
\begin{equation}
\label{4.2}
\frac{\Gamma \left( \alpha \left( k+1 \right)+1 \right)}{\Gamma \left( \alpha k+1 \right)}Y\left( k+1 \right)=2Y\left( k \right)-\sum\limits_{{{k}_{1}}=0}^{k}{Y\left( {{k}_{1}} \right)Y\left( k-{{k}_{1}} \right)}+\delta \left( k \right),
\end{equation}
where $Y\left(0\right)=y\left( 0 \right)$ and  $\delta \left( k \right) =\begin {cases} 1,&\text{if  } k=0\\
                                                                                        0,&\text{otherwise}\end{cases}$.
By using (\ref{4.2}) and recalling (\ref{8}), for $\alpha =0.7$, the approximate solution to the Eq. (\ref{4.1}) up to $O\left( {{t}^{3.5}} \right)$, over the first sub-interval ${{I}_{1}}=[0,0.2]$ is as follows:
\begin{equation}
\label{19}
y\left( t \right)=\text{ }1.10\text{ }{{t}^{0.7}}+1.61\text{ }{{t}^{1.4}}+\text{1}\text{.14 }{{\text{t}}^{2.1}}-0.60\text{ }{{t}^{2.8}}\\
-2.54\text{ }{{t}^{3.5}},\text{   }0\le t\le 0.2.
\end{equation}
The initial condition for the differential equation over the second interval could be evaluated by using Eq. (\ref{19}):
\begin{equation}
y\left( 0.2 \right)=\text{ 0}\text{.55}.
\end{equation}
Now, for the second interval ${{I}_{2}}=[0.2,0.4]$ ,the Eq.(\ref{4.1}) will be in the form of
\begin{equation}
\label{21}
     \begin{cases}
     {{D}_{0.2}^{\alpha }}y\left( t \right)=2y-{{y}^{2}}+1,\text{      }&t>0.2\text{   ,    }0<\alpha \le 1 \\
     y\left( 0.2 \right)=0.55 \\
               \end{cases}.
\end{equation}
Considering $Y\left( 0 \right)=y\left( 0.2 \right)$ and following the same process, the approximate solution to the Eq. (\ref{21}) is
\begin{multline}
$$y(t)=\text{ 0}\text{.55 + }1.98{{\left( t-0.2 \right)}^{0.7}}+1.30{{\left( t-0.2 \right)}^{1.4}}\text{-1}\text{.54}{{\left( t-0.2 \right)}^{2.1}}\\
-3.07{{\left( t-0.2 \right)}^{2.8}}+0.66{{\left( t-0.2 \right)}^{3.5}}\text{, }0.2\le t\le 0.4.$$
\end{multline}
 Figure \ref{fig1} shows a comparison between the result of MSGDTM and the solution to Eq.(\ref{4.1}), with using fractional Adams-Bashforth-Moulton method, for which MATLAB is used to code the algorithm. The numerical solution is also compared with the result obtained by the MATLAB code fde12.m (see \cite{garrappa2011predictor}). The result of MSGDTM and fractional Adams-Bashforth-Moulton method coincide only for the first sub-interval.
As previously mentioned, neglecting the non-local property of the fractional differentiation operators, the MSGDTM results in solving wrong fractional differential equation at the second sub-interval and this, therefore, causes a discontinuity to occur in the first derivative of the solution of fractional differential equation at the beginning of the second sub-interval. It may be worth stating that increasing the number of sub-intervals makes this inaccuracy become larger and in addition, causes discontinuities of the first derivative of the solution to be hidden in the graph (see Figure \ref{fig2}). A comparison of the fractional Adams-Bashforth-Moulton method and the MSGDTM , for the Eq. (\ref{4.1}), is shown in Figure \ref{fig2}, where the interval $I=[0,3]$ has been divided into $M=300$ sub-intervals.
\section{Conclusion}
\label{conclusion}
In this note, it has been briefly mentioned that the solution to the Volterra integral equation (\ref{volterra}) has been definitely proved to have an asymptotic expansion in terms of mixed powers, i.e. integer and non-integer powers whilst under general conditions, the GDTM might result in a power series not involving the integer-order powers. Omitting the integer-order powers from the solution to the fractional differential equation leads to a poor approximation of the exact solution. Moreover, it has been concisely demonstrated that the multi-step generalized differential transform method fails to correctly solve a fractional differential equation due to the fact that the method neglects the effect of non-local structure of fractional differentiation operators. The results have been compared with those from the PECE method of Adams-Bashforth-Moulton type for the solution to the fractional differential equations. It must be mentioned that in the case of $\alpha =1$, the fractional differential equation (\ref{ini1}) reduces to an ordinary differential equation, therefore, the MSGDTM results in the correct solution, which is in a perfect agreement with other classical methods like the fourth order Runge-Kutta method.
\section*{Acknowledgement}
The authors would like to thank Prof. Paola Loreti and Dr. Mirko D'Ovidio for their helpful suggestions.

%\section*{References}

\bibliographystyle{model3-num-names}
\bibliography{mybib}

\begin{thebibliography}{32}
\providecommand{\natexlab}[1]{#1}
\providecommand{\url}[1]{\texttt{#1}}
\providecommand{\href}[2]{#2}
\providecommand{\path}[1]{#1}
\providecommand{\eprint}[1]{\href{http://arxiv.org/abs/#1}{\path{#1}}}
\providecommand{\DOIprefix}{doi:}
\providecommand{\ArXivprefix}{arXiv:}
\providecommand{\URLprefix}{URL: }
\providecommand{\Pubmedprefix}{pmid:}
\providecommand{\doi}[1]{\href{http://dx.doi.org/#1}{\path{#1}}}
\providecommand{\Pubmed}[1]{\href{pmid:#1}{\path{#1}}}
\providecommand{\BIBand}{and}
\providecommand{\bibinfo}[2]{#2}
\ifx\xfnm\undefined \def\xfnm[#1]{\unskip,\space#1}\fi
%Type = Article
\bibitem[{Diethelm and Freed(1998)}]{diethelm1998fracpece}
\bibinfo{author}{Diethelm\xfnm[ K.]}, \bibinfo{author}{Freed\xfnm[ A.D.]}.
\newblock \bibinfo{title}{The frac{PECE} subroutine for the numerical solution
  of differential equations of fractional order}.
\newblock \bibinfo{journal}{Forschung und wissenschaftliches Rechnen}
  \bibinfo{year}{1998};\bibinfo{volume}{1999}:\bibinfo{pages}{57--71}.
%Type = Book
\bibitem[{Carpinteri and Mainardi(2014)}]{carpinteri2014fractals}
\bibinfo{author}{Carpinteri\xfnm[ A.]}, \bibinfo{author}{Mainardi\xfnm[ F.]}.
\newblock \bibinfo{title}{Fractals and fractional calculus in continuum
  mechanics}; vol. \bibinfo{volume}{378}.
\newblock \bibinfo{publisher}{Springer}; \bibinfo{year}{2014}.
%Type = Article
\bibitem[{Atangana and Baleanu(2016)}]{atangana2016new}
\bibinfo{author}{Atangana\xfnm[ A.]}, \bibinfo{author}{Baleanu\xfnm[ D.]}.
\newblock \bibinfo{title}{New fractional derivatives with nonlocal and
  non-singular kernel: theory and application to heat transfer model}.
\newblock \bibinfo{journal}{arXiv preprint arXiv:160203408}
  \bibinfo{year}{2016};.
%Type = Article
\bibitem[{N{\"a}sholm and Holm(2013)}]{nasholm2013fractional}
\bibinfo{author}{N{\"a}sholm\xfnm[ S.P.]}, \bibinfo{author}{Holm\xfnm[ S.]}.
\newblock \bibinfo{title}{On a fractional zener elastic wave equation}.
\newblock \bibinfo{journal}{Fractional Calculus and Applied Analysis}
  \bibinfo{year}{2013};\bibinfo{volume}{16}(\bibinfo{number}{1}):\bibinfo{pages}{26--50}.
%Type = Article
\bibitem[{Lopes et~al.(2013)Lopes, Machado, Pinto and
  Galhano}]{lopes2013fractional}
\bibinfo{author}{Lopes\xfnm[ A.M.]}, \bibinfo{author}{Machado\xfnm[ J.T.]},
  \bibinfo{author}{Pinto\xfnm[ C.M.]}, \bibinfo{author}{Galhano\xfnm[ A.M.]}.
\newblock \bibinfo{title}{Fractional dynamics and {MDS} visualization of
  earthquake phenomena}.
\newblock \bibinfo{journal}{Computers \& Mathematics with Applications}
  \bibinfo{year}{2013};\bibinfo{volume}{66}(\bibinfo{number}{5}):\bibinfo{pages}{647--658}.
%Type = Article
\bibitem[{Koh and Kelly(1990)}]{koh1990application}
\bibinfo{author}{Koh\xfnm[ C.G.]}, \bibinfo{author}{Kelly\xfnm[ J.M.]}.
\newblock \bibinfo{title}{Application of fractional derivatives to seismic
  analysis of base-isolated models}.
\newblock \bibinfo{journal}{Earthquake engineering \& structural dynamics}
  \bibinfo{year}{1990};\bibinfo{volume}{19}(\bibinfo{number}{2}):\bibinfo{pages}{229--241}.
%Type = Article
\bibitem[{Rossikhin and Shitikova(1997)}]{rossikhin1997applications}
\bibinfo{author}{Rossikhin\xfnm[ Y.A.]}, \bibinfo{author}{Shitikova\xfnm[
  M.V.]}.
\newblock \bibinfo{title}{Applications of fractional calculus to dynamic
  problems of linear and nonlinear hereditary mechanics of solids}.
\newblock \bibinfo{journal}{Applied Mechanics Reviews}
  \bibinfo{year}{1997};\bibinfo{volume}{50}:\bibinfo{pages}{15--67}.
%Type = Article
\bibitem[{Rossikhin and Shitikova(2010)}]{rossikhin2010application}
\bibinfo{author}{Rossikhin\xfnm[ Y.A.]}, \bibinfo{author}{Shitikova\xfnm[
  M.V.]}.
\newblock \bibinfo{title}{Application of fractional calculus for dynamic
  problems of solid mechanics: novel trends and recent results}.
\newblock \bibinfo{journal}{Applied Mechanics Reviews}
  \bibinfo{year}{2010};\bibinfo{volume}{63}(\bibinfo{number}{1}):\bibinfo{pages}{010801}.
%Type = Article
\bibitem[{Agrawal(2004)}]{agrawal2004general}
\bibinfo{author}{Agrawal\xfnm[ O.P.]}.
\newblock \bibinfo{title}{A general formulation and solution scheme for
  fractional optimal control problems}.
\newblock \bibinfo{journal}{Nonlinear Dynamics}
  \bibinfo{year}{2004};\bibinfo{volume}{38}(\bibinfo{number}{1-4}):\bibinfo{pages}{323--337}.
%Type = Article
\bibitem[{Jelicic and Petrovacki(2009)}]{jelicic2009optimality}
\bibinfo{author}{Jelicic\xfnm[ Z.D.]}, \bibinfo{author}{Petrovacki\xfnm[ N.]}.
\newblock \bibinfo{title}{Optimality conditions and a solution scheme for
  fractional optimal control problems}.
\newblock \bibinfo{journal}{Structural and Multidisciplinary Optimization}
  \bibinfo{year}{2009};\bibinfo{volume}{38}(\bibinfo{number}{6}):\bibinfo{pages}{571--581}.
%Type = Article
\bibitem[{Magin(2010)}]{magin2010fractional}
\bibinfo{author}{Magin\xfnm[ R.L.]}.
\newblock \bibinfo{title}{Fractional calculus models of complex dynamics in
  biological tissues}.
\newblock \bibinfo{journal}{Computers \& Mathematics with Applications}
  \bibinfo{year}{2010};\bibinfo{volume}{59}(\bibinfo{number}{5}):\bibinfo{pages}{1586--1593}.
%Type = Article
\bibitem[{Ding and Ye(2009)}]{ding2009fractional}
\bibinfo{author}{Ding\xfnm[ Y.]}, \bibinfo{author}{Ye\xfnm[ H.]}.
\newblock \bibinfo{title}{A fractional-order differential equation model of
  {HIV} infection of {CD4}+ {T}-cells}.
\newblock \bibinfo{journal}{Mathematical and Computer Modelling}
  \bibinfo{year}{2009};\bibinfo{volume}{50}(\bibinfo{number}{3}):\bibinfo{pages}{386--392}.
%Type = Article
\bibitem[{Gl{\"o}ckle and Nonnenmacher(1995)}]{glockle1995fractional}
\bibinfo{author}{Gl{\"o}ckle\xfnm[ W.G.]}, \bibinfo{author}{Nonnenmacher\xfnm[
  T.F.]}.
\newblock \bibinfo{title}{A fractional calculus approach to self-similar
  protein dynamics}.
\newblock \bibinfo{journal}{Biophysical Journal}
  \bibinfo{year}{1995};\bibinfo{volume}{68}(\bibinfo{number}{1}):\bibinfo{pages}{46--53}.
%Type = Book
\bibitem[{Magin(2006)}]{magin2006fractional}
\bibinfo{author}{Magin\xfnm[ R.L.]}.
\newblock \bibinfo{title}{Fractional calculus in bioengineering}.
\newblock \bibinfo{publisher}{Begell House Redding}; \bibinfo{year}{2006}.
%Type = Article
\bibitem[{Diethelm and Walz(1997)}]{diethelm1997numerical}
\bibinfo{author}{Diethelm\xfnm[ K.]}, \bibinfo{author}{Walz\xfnm[ G.]}.
\newblock \bibinfo{title}{Numerical solution of fractional order differential
  equations by extrapolation}.
\newblock \bibinfo{journal}{Numerical algorithms}
  \bibinfo{year}{1997};\bibinfo{volume}{16}(\bibinfo{number}{3}):\bibinfo{pages}{231--253}.
%Type = Article
\bibitem[{Chen and Liu(2008)}]{chen2008adi}
\bibinfo{author}{Chen\xfnm[ S.]}, \bibinfo{author}{Liu\xfnm[ F.]}.
\newblock \bibinfo{title}{{ADI}-euler and extrapolation methods for the
  two-dimensional fractional advection-dispersion equation}.
\newblock \bibinfo{journal}{Journal of Applied Mathematics and Computing}
  \bibinfo{year}{2008};\bibinfo{volume}{26}(\bibinfo{number}{1-2}):\bibinfo{pages}{295--311}.
%Type = Article
\bibitem[{Lubich(1985)}]{lubich1985fractional}
\bibinfo{author}{Lubich\xfnm[ C.]}.
\newblock \bibinfo{title}{Fractional linear multistep methods for
  {A}bel-{V}olterra integral equations of the second kind}.
\newblock \bibinfo{journal}{Mathematics of computation}
  \bibinfo{year}{1985};\bibinfo{volume}{45}(\bibinfo{number}{172}):\bibinfo{pages}{463--469}.
%Type = Article
\bibitem[{Galeone and Garrappa(2006)}]{galeone2006multistep}
\bibinfo{author}{Galeone\xfnm[ L.]}, \bibinfo{author}{Garrappa\xfnm[ R.]}.
\newblock \bibinfo{title}{On multistep methods for differential equations of
  fractional order}.
\newblock \bibinfo{journal}{Mediterranean Journal of Mathematics}
  \bibinfo{year}{2006};\bibinfo{volume}{3}(\bibinfo{number}{3}):\bibinfo{pages}{565--580}.
%Type = Article
\bibitem[{Galeone and Garrappa(2009)}]{galeone2009explicit}
\bibinfo{author}{Galeone\xfnm[ L.]}, \bibinfo{author}{Garrappa\xfnm[ R.]}.
\newblock \bibinfo{title}{Explicit methods for fractional differential
  equations and their stability properties}.
\newblock \bibinfo{journal}{Journal of Computational and Applied Mathematics}
  \bibinfo{year}{2009};\bibinfo{volume}{228}(\bibinfo{number}{2}):\bibinfo{pages}{548--560}.
%Type = Book
\bibitem[{Baleanu et~al.(2016)Baleanu, Diethelm, Scalas and
  Trujillo}]{baleanu2016fractional}
\bibinfo{author}{Baleanu\xfnm[ D.]}, \bibinfo{author}{Diethelm\xfnm[ K.]},
  \bibinfo{author}{Scalas\xfnm[ E.]}, \bibinfo{author}{Trujillo\xfnm[ J.J.]}.
\newblock \bibinfo{title}{Fractional calculus: models and numerical methods};
  vol.~\bibinfo{volume}{5}.
\newblock \bibinfo{publisher}{World Scientific}; \bibinfo{year}{2016}.
%Type = Article
\bibitem[{Diethelm(1997)}]{diethelm1997algorithm}
\bibinfo{author}{Diethelm\xfnm[ K.]}.
\newblock \bibinfo{title}{An algorithm for the numerical solution of
  differential equations of fractional order}.
\newblock \bibinfo{journal}{Electronic transactions on numerical analysis}
  \bibinfo{year}{1997};\bibinfo{volume}{5}(\bibinfo{number}{1}):\bibinfo{pages}{1--6}.
%Type = Article
\bibitem[{Diethelm et~al.(2004)Diethelm, Ford and Freed}]{diethelm2004detailed}
\bibinfo{author}{Diethelm\xfnm[ K.]}, \bibinfo{author}{Ford\xfnm[ N.J.]},
  \bibinfo{author}{Freed\xfnm[ A.D.]}.
\newblock \bibinfo{title}{Detailed error analysis for a fractional adams
  method}.
\newblock \bibinfo{journal}{Numerical algorithms}
  \bibinfo{year}{2004};\bibinfo{volume}{36}(\bibinfo{number}{1}):\bibinfo{pages}{31--52}.
%Type = Article
\bibitem[{Diethelm(2003)}]{diethelm2003efficient}
\bibinfo{author}{Diethelm\xfnm[ K.]}.
\newblock \bibinfo{title}{Efficient solution of multi-term fractional
  differential equations using {P} ({EC}) m{E} methods}.
\newblock \bibinfo{journal}{Computing}
  \bibinfo{year}{2003};\bibinfo{volume}{71}(\bibinfo{number}{4}):\bibinfo{pages}{305--319}.
%Type = Article
\bibitem[{Garrappa(2011)}]{garrappa2011predictor}
\bibinfo{author}{Garrappa\xfnm[ R.]}.
\newblock \bibinfo{title}{Predictor-corrector {PECE} method for fractional
  differential equations}.
\newblock \bibinfo{journal}{MATLAB Central File Exchange [File ID: 32918]}
  \bibinfo{year}{2011};.
%Type = Book
\bibitem[{Diethelm(2010)}]{diethelm2010analysis}
\bibinfo{author}{Diethelm\xfnm[ K.]}.
\newblock \bibinfo{title}{The analysis of fractional differential equations: An
  application-oriented exposition using differential operators of Caputo type}.
\newblock \bibinfo{publisher}{Springer}; \bibinfo{year}{2010}.
%Type = Article
\bibitem[{Arikoglu and Ozkol(2007)}]{arikoglu2007solution}
\bibinfo{author}{Arikoglu\xfnm[ A.]}, \bibinfo{author}{Ozkol\xfnm[ I.]}.
\newblock \bibinfo{title}{Solution of fractional differential equations by
  using differential transform method}.
\newblock \bibinfo{journal}{Chaos, Solitons \& Fractals}
  \bibinfo{year}{2007};\bibinfo{volume}{34}(\bibinfo{number}{5}):\bibinfo{pages}{1473--1481}.
%Type = Article
\bibitem[{Odibat et~al.(2008)Odibat, Momani and Erturk}]{odibat2008generalized}
\bibinfo{author}{Odibat\xfnm[ Z.]}, \bibinfo{author}{Momani\xfnm[ S.]},
  \bibinfo{author}{Erturk\xfnm[ V.S.]}.
\newblock \bibinfo{title}{Generalized differential transform method:
  application to differential equations of fractional order}.
\newblock \bibinfo{journal}{Applied Mathematics and Computation}
  \bibinfo{year}{2008};\bibinfo{volume}{197}(\bibinfo{number}{2}):\bibinfo{pages}{467--477}.
%Type = Article
\bibitem[{Kilbas et~al.(2006)Kilbas, Srivastava and
  Trujillo}]{kilbas2006theory}
\bibinfo{author}{Kilbas\xfnm[ A.]}, \bibinfo{author}{Srivastava\xfnm[ H.]},
  \bibinfo{author}{Trujillo\xfnm[ J.]}.
\newblock \bibinfo{title}{Theory and applications of fractional differential
  equations, volume 204 (north-holland mathematics studies)}
  \bibinfo{year}{2006};.
%Type = Article
\bibitem[{Odibat(2008)}]{odibat2008differential}
\bibinfo{author}{Odibat\xfnm[ Z.M.]}.
\newblock \bibinfo{title}{Differential transform method for solving volterra
  integral equation with separable kernels}.
\newblock \bibinfo{journal}{Mathematical and Computer Modelling}
  \bibinfo{year}{2008};\bibinfo{volume}{48}(\bibinfo{number}{7}):\bibinfo{pages}{1144--1149}.
%Type = Article
\bibitem[{Bervillier(2012)}]{bervillier2012status}
\bibinfo{author}{Bervillier\xfnm[ C.]}.
\newblock \bibinfo{title}{Status of the differential transformation method}.
\newblock \bibinfo{journal}{Applied Mathematics and Computation}
  \bibinfo{year}{2012};\bibinfo{volume}{218}(\bibinfo{number}{20}):\bibinfo{pages}{10158--10170}.
%Type = Article
\bibitem[{Lubich(1983)}]{lubich1983runge}
\bibinfo{author}{Lubich\xfnm[ C.]}.
\newblock \bibinfo{title}{Runge-{Kutta} theory for {Volterra} and {Abel}
  integral equations of the second kind}.
\newblock \bibinfo{journal}{Mathematics of computation}
  \bibinfo{year}{1983};\bibinfo{volume}{41}(\bibinfo{number}{163}):\bibinfo{pages}{87--102}.
%Type = Article
\bibitem[{Ert{\"u}rk et~al.(2011)Ert{\"u}rk, Odibat and
  Momani}]{erturk2011approximate}
\bibinfo{author}{Ert{\"u}rk\xfnm[ V.S.]}, \bibinfo{author}{Odibat\xfnm[ Z.M.]},
  \bibinfo{author}{Momani\xfnm[ S.]}.
\newblock \bibinfo{title}{An approximate solution of a fractional order
  differential equation model of human {T-cell} lymphotropic virus {I}
  ({HTLV-I}) infection of {T}-cells}.
\newblock \bibinfo{journal}{Computers \& Mathematics with applications}
  \bibinfo{year}{2011};\bibinfo{volume}{62}(\bibinfo{number}{3}):\bibinfo{pages}{996--1002}.

\end{thebibliography}

%\endinput

\hfill \break

   Department of Structure and Material\\
Universiti Teknologi Malaysia\\
81310 Skudai, Johor, Malaysia\\[4pt]
e-mail: alireza.momenzadeh@gmail.com
   \\[12pt]

    $^*$ Dipartimento di Scienze di Base e Applicate per l'Ingegneria \\
Sapienza Universit\`{a} di Roma \\
Via Antonio Scarpa n. 16 \\
00161 Rome, Italy \\[4pt]
  e-mail: sima.sarvahrabi@sbai.uniroma1.it
\end{document}